\documentclass[11pt]{article}

\usepackage{amssymb}

\newcommand{\newc}{\newcommand}

\newc{\eqnoset}{\setcounter{equation}{0}}

\newcommand{\mref}[1]{(\ref{#1})}
\newcommand{\reflemm}[1]{Lemma~\ref{#1}}
\newcommand{\refrem}[1]{Remark~\ref{#1}}
\newcommand{\reftheo}[1]{Theorem~\ref{#1}}

\newcommand{\refcoro}[1]{Corollary~\ref{#1}}

\newcommand{\refsec}[1]{Section~\ref{#1}}

\newcommand{\beq}{\begin{equation}}
\newcommand{\eeq}{\end{equation}}
\newcommand{\beqno}[1]{\begin{equation}\label{#1}}

\newcommand{\barr}{\begin{array}}
\newcommand{\earr}{\end{array}}

\newc{\bearr}{\begin{eqnarray*}}
\newc{\eearr}{\end{eqnarray*}}

\newc{\bearrno}[1]{\begin{eqnarray}\label{#1}}
\newc{\eearrno}{\end{eqnarray}}

\newc{\non}{\nonumber}
\newc{\nol}{\nonumber\nl}

\newcommand{\bdes}{\begin{description}}
\newcommand{\edes}{\end{description}}
\newc{\benu}{\begin{enumerate}}
\newc{\eenu}{\end{enumerate}}
\newc{\btab}{\begin{tabular}}
\newc{\etab}{\end{tabular}}

\newtheorem{theorem}{Theorem}[section]
\newtheorem{defi}[theorem]{Definition}
\newtheorem{lemma}[theorem]{Lemma}
\newtheorem{rem}[theorem]{Remark}
\newtheorem{exam}[theorem]{Example}
\newtheorem{propo}[theorem]{Proposition}
\newtheorem{corol}[theorem]{Corollary}

\newcommand{\btheo}[1]{\begin{theorem}\label{#1}}
\newc{\brem}[1]{\begin{rem}\label{#1}\em}
\newc{\bexam}[1]{\begin{exam}\label{#1}\em}
\newc{\bdefi}[1]{\begin{defi}\label{#1}}
\newcommand{\blemm}[1]{\begin{lemma}\label{#1}}
\newcommand{\bprop}[1]{\begin{propo}\label{#1}}
\newcommand{\bcoro}[1]{\begin{corol}\label{#1}}
\newcommand{\etheo}{\end{theorem}}
\newcommand{\elemm}{\end{lemma}}
\newcommand{\eprop}{\end{propo}}
\newcommand{\ecoro}{\end{corol}}
\newc{\erem}{\end{rem}}
\newc{\eexam}{\end{exam}}
\newc{\edefi}{\end{defi}}

\newc{\rmk}[1]{{\bf REMARK #1: }}
\newc{\DN}[1]{{\bf DEFINITION #1: }}

\newcommand{\bproof}{{\bf Proof:~~}}
\newc{\eproof}{{\vrule height8pt width5pt depth0pt}\vspace{3mm}}

\newc{\bfrac}[2]{\dspl{\frac{#1}{#2}}}

\newc{\nid}{\noindent}

\newcommand{\dspl}{\displaystyle}
\newc{\grad}{\nabla}
\newc{\Div}{\mbox{div}}
\newc{\pdt}[1]{\dspl{\frac{\partial{#1}}{\partial t}}}
\newc{\pdn}[1]{\dspl{\frac{\partial{#1}}{\partial \nu}}}
\newc{\pdNi}[1]{\dspl{\frac{\partial{#1}}{\partial \mathcal{N}_i}}}
\newc{\pD}[2]{\dspl{\frac{\partial{#1}}{\partial #2}}}
\newc{\dt}{\dspl{\frac{d}{dt}}}
\newc{\bdry}[1]{\mbox{$\partial #1$}}
\newc{\sgn}{\mbox{sign}}

\newc{\Hess}[1]{\frac{\partial^2 #1}{\pdh z_i \pdh z_j}}
\newc{\hess}[1]{\partial^2 #1/\pdh z_i \pdh z_j}

\newc{\ag}{\alpha}
\newc{\bg}{\beta}
\newc{\cg}{\gamma}\newc{\Cg}{\Gamma}
\newc{\dg}{\delta}\newc{\Dg}{\Delta}
\newc{\eg}{\varepsilon}
\newc{\zg}{\zeta}
\newc{\thg}{\theta}
\newc{\llg}{\lambda}\newc{\LLg}{\Lambda}
\newc{\kg}{\kappa}
\newc{\rg}{\rho}
\newc{\sg}{\sigma}\newc{\Sg}{\Sigma}
\newc{\tg}{\tau}
\newc{\fg}{\phi}\newc{\Fg}{\Phi}
\newc{\vfg}{\varphi}
\newc{\og}{\omega}\newc{\Og}{\Omega}
\newc{\pdh}{\partial}

\newc{\ccG}{{\cal G}}

\newc{\ii}[1]{\int_{#1}}
\newc{\iidx}[2]{{\dspl\int_{#1}~#2~dx}}
\newc{\bii}[1]{{\dspl \ii{#1} }}
\newc{\biii}[2]{{\dspl \iii{#1}{#2} }}
\newc{\su}[2]{\sum_{#1}^{#2}}
\newc{\bsu}[2]{{\dspl \su{#1}{#2} }}

\newc{\biiom}[1]{{\dspl\int_{\bdrom}~ #1 ~d\sg}}
\newc{\io}[1]{{\dspl\int_{\Og}~ #1 ~dx}}
\newc{\bio}[1]{{\dspl\int_{\bdrom}~ #1 ~d\sg}}
\newc{\bsir}{\bsu{i=1}{r}}
\newc{\bsim}{\bsu{i=1}{m}}

\newc{\iibr}[2]{\iidx{\bprw{#1}}{#2}}
\newc{\Intbr}[1]{\iibr{R}{#1}}
\newc{\intbr}[1]{\iibr{\rg}{#1}}
\newc{\intt}[3]{\int_{#1}^{#2}\int_\Og~#3~dxdt}

\newc{\itQ}[2]{\dspl{\int\hspace{-2.5mm}\int_{#1}~#2~dz}}
\newc{\mitQ}[2]{\dspl{\rule[1mm]{4mm}{.3mm}\hspace{-5.3mm}\int\hspace{-2.5mm}\int_{#1}~#2~dz}}
\newc{\mitQQ}[3]{\dspl{\rule[1mm]{4mm}{.3mm}\hspace{-5.3mm}\int\hspace{-2.5mm}\int_{#1}~#2~#3}}

\newc{\mitx}[2]{\dspl{\rule[1mm]{3mm}{.3mm}\hspace{-4mm}\int_{#1}~#2~dx}}
\newc{\mitmu}[2]{\dspl{\rule[1mm]{3mm}{.3mm}\hspace{-4mm}\int_{#1}~#2~d\mu}}
\newc{\iidmu}[2]{{\dspl\int_{#1}~#2~d\mu}}

\newc{\iidm}[3]{{\dspl\int_{#1}~#2~d #3}}

\newc{\itQmu}[2]{\dspl{\int\hspace{-2.5mm}\int_{#1}~#2~d\mu}}
\newc{\mitQmu}[2]{\dspl{\rule[1mm]{4mm}{.3mm}\hspace{-5.3mm}\int\hspace{-2.5mm}\int_{#1}~#2~d\mu}}

\newc{\mitQq}[2]{\dspl{\rule[1mm]{4mm}{.3mm}\hspace{-5.3mm}\int\hspace{-2.5mm}\int_{#1}~#2~d\bar{z}}}
\newc{\itQq}[2]{\dspl{\int\hspace{-2.5mm}\int_{#1}~#2~d\bar{z}}}

\newc{\pder}[2]{\dspl{\frac{\partial #1}{\partial #2}}}

\newc{\bdrom}{\bdry{\Og}}

\newc{\bilhom}{\mbox{Bil}(\mbox{Hom}(\RR^{nm},\RR^{nm}))}
\newc{\VV}[1]{{V(Q_{#1})}}

\newc{\ccA}{{\mathcal A}}
\newc{\ccB}{{\mathcal B}}
\newc{\ccC}{{\mathcal C}}
\newc{\ccD}{{\mathcal D}}
\newc{\ccE}{{\mathcal E}}
\newc{\ccH}{\mathcal{H}}
\newc{\ccF}{\mathcal{F}}
\newc{\ccI}{{\mathcal I}}
\newc{\ccJ}{{\mathcal J}}
\newc{\ccK}{{\mathcal K}}
\newc{\ccP}{{\mathcal P}}
\newc{\ccQ}{{\mathcal Q}}
\newc{\ccR}{{\mathcal R}}
\newc{\ccS}{{\mathcal S}}
\newc{\ccT}{{\mathcal T}}
\newc{\ccX}{{\mathcal X}}
\newc{\ccY}{{\mathcal Y}}
\newc{\ccZ}{{\mathcal Z}}

\newc{\bb}[1]{{\mathbf #1}}

\newc{\myprod}[1]{\langle #1 \rangle}
\newc{\mypar}[1]{\left( #1 \right)}

\newc{\BLLg}{\mathbf{\LLg}}

\newc{\mA}{\mathbf{A}}
\newc{\mB}{\mathbf{B}}
\newc{\mC}{\mathbf{C}}
\newc{\mD}{\mathbf{D}}
\newc{\mE}{\mathbf{E}}
\newc{\mF}{\mathbf{F}}
\newc{\mJ}{\mathbf{J}}
\newc{\mG}{\mathbf{G}}
\newc{\mP}{\mathbf{P}}
\newc{\mR}{\mathbf{R}}
\newc{\mQ}{\mathbf{Q}}
\newc{\mX}{\mathbf{X}}
\newc{\muu}{\mathbf{u}}
\newc{\mvv}{\mathbf{v}}

\newc{\mllg}{\mathbb{\lambda}}
\newc{\mLLg}{\mathbf{\LLg}}

\newc{\lspn}[2]{\mbox{$\| #1\|_{\Lsp{#2}}$}}
\newc{\Lpn}[2]{\mbox{$\| #1\|_{#2}$}}
\newc{\Hn}[1]{\mbox{$\| #1\|_{H^1(\Og)}$}}

\newc{\mynorm}[2]{\| #1\|_{#2}}

\newcommand{\RR}{{\rm I\kern -1.6pt{\rm R}}}

\newc{\itQQ}[2]{\dspl{\int_{#1}#2\,dz}}
\newc{\mmitQQ}[2]{\dspl{\rule[1mm]{4mm}{.3mm}\hspace{-4.3mm}\int_{#1}~#2~dz}}
\newc{\MmitQQ}[2]{\dspl{\rule[1mm]{4mm}{.3mm}\hspace{-4.3mm}\int_{#1}~#2~d\mu}}

\newc{\MUmitQQ}[3]{\dspl{\rule[1mm]{4mm}{.3mm}\hspace{-4.3mm}\int_{#1}~#2~d#3}}
\newc{\MUitQQ}[3]{\dspl{\int_{#1}~#2~d#3}}

\newc{\mccP}{\mathbb{P}}
\newc{\mccK}{\mathbb{K}}

\newc{\DKTmU}{\mccK(U)}
\newc{\DKTmUold}{(K_U(U)^{-1})^T}

\newc{\myPi}{\mathbf{W}}
\newc{\myIbar}{\bar{\ccI}_1}
\newc{\myIhat}{\hat{\ccI}_1}
\newc{\myIbreve}{\breve{\ccI}_0}

\newc{\mmk}{\mathbf{k}}

\newcommand{\ma}{\mathbf{a}}
\newcommand{\mg}{\mathbf{g}}

\newc{\mfu}{\mathbf{f_u}}
\newc{\mh}{\mathbf{h}}
\newc{\mb}{\mathbf{b}}

\newcommand{\barrl}[2]{\barr{ll}\lefteqn{#1}\hspace{#2}&\\}

\newc{\mN}{\mathbf{N}}
\newc{\mI}{\mathbf{I}}
\newc{\mH}{\mathbf{H}}

\newc{\mk}{\mathbf{k}}
\newc{\mr}{\mathbf{r}}

\newc{\DIAGM}[2]{\left[\barr{ccc}#1&0\ldots&0\\
	\vdots&\ddots&\vdots\\0&\ldots0&#2\earr \right]}
\newc{\DiagM}[2]{\mbox{diag}\left[#1
	\cdots #2 \right]}
\newc{\vVEC}[2]{\left[\barr{c}#1\\
	\vdots\\#2\earr \right]}
\newc{\hVEC}[2]{\left[#1
	\cdots #2 \right]}

\newc{\mq}{\mathbf{q}}

\newc{\msys}[1]{\left\{\barr{l}#1\earr
	\right.}
\newc{\msysa}[1]{\left\{\barr{ll}#1\earr
	\right.}
\newc{\mL}{\mathbb{L}}

\newc{\bbM}{\mathbb{M}}
\newc{\mat}[1]{\left[\barr{cc}#1\earr\right]}

\begin{document}

\vspace*{-.8in}
\begin{center} {\LARGE\em On the Global Existence of a Class of Strongly Coupled Parabolic Systems.}

 \end{center}

\vspace{.1in}

\begin{center}

{\sc Dung Le}{\footnote {Department of Mathematics, University of
Texas at San
Antonio, One UTSA Circle, San Antonio, TX 78249. {\tt Email: Dung.Le@utsa.edu}
}}

\end{center}

\begin{abstract}
We establish the  existence of strong solutions to a class of cross diffusion systems on $\RR^N$ consists of $m$ equations ($m,N\ge 2$). which generalizes the Shigesada-Kawasaki-Teramoto (SKT) model in population dynamics. We introduce the concept of a {\em strong-weak solution} of the systems and show that their existence can be established under weaker conditions. These {\em strong-weak solutions} coincide with strong solutions so that the existence of strong solutions is proved. The SKT model on planar domains ($N=2$) with cubic diffusions and advections is completely solved.  \end{abstract}

\vspace{.2in}

\section{Introduction} \label{intro}\eqnoset

In this paper, let $\Og$ be a bounded domain in $\RR^N$, $N\ge2$, with smooth boundary $\partial \Og$, and $T>0$.  We study the solvability of the strongly coupled parabolic system
\beqno{fullsys0a}\msysa{W_t=\Div(\ma(W) DW)+\Div(\hat{\mb}(W))+\mb(W) DW+\mg(W) &\mbox{in $Q=\Og\times(0,T)$,}\\ \mbox{Homogeneous Dirichlet or Neumann boundary conditions}&\mbox{on $\partial\Og\times(0,T)$,}\\W=u_0&\mbox{on $\Og\times(0,T)$.}}\eeq
Here,  $W=[u_i]_{i=1}^m$, a vector in $\RR^m$ and $\ma,\mb,\mg$ are square matrices of size $m\times m$ for some $m\ge2$. The entries of these matrices are functions in $W$. The initial data  $u_0$ is a given vector valued function in $W^{1,N_0}(\Og)$ for some $N_0>N$. For simplicity we will assume the components of $\hat{\mb}(W),\mb(W),\mg(W)$ have polynomial growths in $W$ throughout this paper although many results in this paper holds under the assumption that they are bounded if $W$ is bounded.

According to the usual definition, a {\em strong} solution to \mref{fullsys0a} is a vector valued function $W\in W^{2,2}_{loc}(Q)$ which has bounded derivative $DW$ and solves \mref{fullsys0a} a.e. in $Q$.
A {\em weak} solution to \mref{fullsys0a} is a vector valued function $W\in V_2(Q)$ (see \cite{LSU}) which satisfies the integral form of \mref{fullsys0a}. That is, for all $\psi\in C^1(Q)$ and any $Q_t=\Og\times(0,t)$ with $t\in(0,T)$
$$\left.\iidx{\Og}{W\psi}\right|_{t=0}^{t}-\itQ{Q_t}{W\psi_t}+\itQ{Q_t}{\ma DW D\psi}=-\itQ{Q_t}{(\hat{\mb} D\psi+\mb DW\psi +\mg \psi)}.$$

The existence problem of a (unique) strong solution to \mref{fullsys0a} was investigated by Amann. He uses interpolation functional space theory and shows that if the parameters of the regular parabolic \mref{fullsys0a} are bounded and  \beqno{est1a} \sup_{(0,T)}\|W\|_{W^{1,2p}(\Og)} <M \mbox{ for some $p>N/2$ and $M$}\eeq then there is a unique strong solution of \mref{fullsys0a}. The existence problem of a weak solution to \mref{fullsys0a} was proved easily by (for example) appropriate Gal\"erkin methods in literature. However, the uniqueness problem was largely open if the coefficients of the system depend on $W$.

For nonlinear strongly coupled systems like \mref{fullsys0a} one would start by proving the boundedness of solutions because Amann's theory worked with bounded $\mA$. For nonlinear strongly coupled systems like \mref{fullsys0a} this problem would be a very hard one already. The next obstacle is, and even harder,  the estimate of higher order norms like \mref{est1a}. This problem is closely related to the regularity of parabolic systems.

In this paper, we introduce the concept of a {\em strong-weak solution} to \mref{fullsys0a}. We say that a vector valued function $W$ is a strong weak solution of \mref{fullsys0a} if  $W\in W^{1,p}(Q)$ (whose spatial derivative $DW\in L^p(Q)$) for some $p>N$ and solves \mref{fullsys0a} weakly. i.e. for all $\psi\in C^1_0(Q)$
\beqno{wf}-\itQ{Q}{W_t\psi}+\itQ{Q}{\ma DW D\psi}=-\itQ{Q}{\hat{\mb} D\psi}+\itQ{Q}{(\mb DW +\mg )\psi}.\eeq  The temporal derivative $W_t$ could be replaced by the Steklov average of $W$.

Of course, the concept of a strong weak solution is weaker than that of strong solution and stronger than that of weak solutions. The advantage of this definition is that a strong weak solution is unique if it exists.

In this paper, we will establish the existence of a unique strong weak solution under an integral assumption which is weaker than \mref{est1a} but requires some extra structure of conditions on the systems. In particular, the spectral gap condition which requires the eigenvalues of $\ma$ are not too far apart.
In this case, we see that the condition \mref{est1a} of Amann can be replaced by a weaker one
\beqno{est1aa}\itQ{\Og\times(0,T)}{|DW|^{2p}}<M \mbox{ for some $p>N/2$ and $M$}\eeq for {\em all} strong weak solutions of \mref{fullsys0a}. In the process of establishing this, we can assume that strong weak solution are H\"older continuous. We can still assert that the {\em strong} solution exists on $\Og\times(0,T)$.

We always assume that $\ma$ satisfies the ellipticity condition (this comes from the normal ellipticity of Amann). That is,  there are some function $\llg(W)$ and a positive constant $\llg_0$ such that $\llg(W)\ge\llg_0$ and
\beqno{ellcond}\myprod{\ma DW,DW}\ge \llg(W)|DW|^2 \quad \forall W\in C^1(\Og,\RR^{mN}).\eeq

Let $\llg\le \LLg$ be smallest and largest eigenvalues of a square matrix $\ma$. We look at the ratio $\nu_*=\frac{\llg}{\LLg}$. We proved in \cite{letrans, dleANS} that if $s>-1$ and $\frac{s}{s+2}<\nu_*$ then \mref{ellcond} implies a $c>0$ such that
\beqno{ellcond10}  \myprod{\ma DX,D(|X|^{s}X)}\ge c\llg|X|^{s}|DX|^2, \quad \forall X\in C^1(\Og,\RR^{mN}).\eeq

Accordingly, we will say that the system satisfies {\em a spectral gap condition for some $p>1/2$}  if the above holds for  the matrix $\ma$ and $s=2p-2$. That is $1-1/p <\nu_*$. In this case, \mref{ellcond10} becomes
$$\myprod{\ma DX,D(|X|^{2p-2}X)}\ge c\llg|X|^{2p-2}|DX|^2, \quad \forall X\in C^1(\Og,\RR^{mN}).$$

Of course, if $\nu_*>1-2/N$ then the spectral gap condition holds for some $p>N/2$.

We combine with Amann's results. For $u_0\in W^{1,p}(\Og)$ we can find a unique strong solution $W_*$ which  exists in $(0,\dg_0)$ for some $\dg_0>0$. We then set up fixed point method in some appropriate space $X$, which will be defined later.
Denote $W_0=W_*(x,\dg_0)$. We define the map $\mL(\Fg)=W$ where $W$ the unique weak solution solution of the {\em linear} problem
\beqno{Ldef}\msysa{W_t=\Div(\ma(\Fg)DW)+\Div(\hat{\mb}(\Fg))+\mb(\Fg)DW+\mg(\Fg)& \mbox{ in }\Og\times(\dg_0,T),\\
W(x,\dg_0)=W_0(x)&\mbox{ on }\Og,\\
\mbox{Homogeneous Dirichlet or Neumann conditions}&\mbox{ on }\partial\Og\times (\dg_0,T).}
\eeq

By the well known results of Ladyzhenskaya {\em et al} in \cite{LSU} it is well known that this has a unique weak solution $W$. The coefficients of the system are smooth in $(0,\dg_0)$ so that $W$ is classical. We need to show that $\mL$ has a unique fixed point in $X$. Gluing this fixed point with the strong solution $W_*$ in $(0,\dg_0)$ we obtain the desired strong weak solution of \mref{fullsys0a}. 
Note that the so glued solution satisfies \mref{wf}.

According to the well known Leray-Schauder theory, we have to consider solutions of the equation $\tau\mL(\Fg)=\Fg$ which is equivalent to  (writing $w=\Fg$)
\beqno{famsys}\msysa{w_t=\Div(\ma(w)Dw)+\tau\Div(\hat{\mb}(w))+\mb(w) Dw+\tau\mg(w) & \mbox{ in $\Og\times(\dg_0,T),$}\\ w(x,0)=\tau W_0(\dg_0,x)&\mbox{ on $\Og$}.}\eeq

In order to apply the Leray-Schauder theorem, we will prove that such $w$ is bounded uniformly for $\tau\in[0,1]$.  There is a technical subtlety in our argument below if we consider \mref{famsys} alone. The estimates, via our parabolic techniques using cutoff functions in $t$, are verified only away from $\dg_0$. To remedy this, we extend $w$ to a solution of system defined on $(0,T)$. The estimates for $w$ on $(\dg_0,T)$ will then be those of this extension away from $0$.

The solution $w$ can be extended backward by $\tau W_*$ in $(0,\dg_0)$ (thanks to the initial condition in \mref{famsys}, $W_*$ is the strong solution, or fix point of \mref{Ldef} in $(0,\dg_0]$) and $w=\tau W_*$ satisfies (a linear system)
\beqno{famsysz}\msysa{w_t=\Div(\ma(W_*)Dw)+\tau\Div(\hat{\mb}(W_*))+\mb(W_*) Dw+\tau\mg(W_*) & \mbox{ in $\Og\times(0,\dg_0),$}\\ w(x,0)=\tau W_*(0,x)=\tau u_0(x)&\mbox{ on $\Og$}.}\eeq

We now consider the gluing solution which solves \mref{famsys} and \mref{famsysz}  in $(0,\dg_0)$ and $(\dg_0,T)$ respectively (so that the gluing $w$ now solves \mref{famsysz} in $(0,T)$). 
Obviously $w$ is smooth in $(0,\dg_0)$ so that we can estimate $w$ in $(\dg_0,T)$ (away from $t=0$) by using cutoff functions in $t$. The coefficients of the system are now nonsmooth and we will be only concerned with $t$ away from $0$.

Thus, the main problems are: \bdes\item[1)] Define a space $X$ such that $\mL:X\to X$  is $\mL$ a compact map. \item[2)] Uniform estimates of solutions to  \mref{famsys} in $X$. \edes

These are the main tasks of this paper and we have the following main result on the existence of a unique strong weak solution in $\Og\times(0,T)$.

\btheo{main} Assume that there is $q>N/2$ such that $1-1/q<\nu_*$ and a constant $M$ such that any strong weak solution $W$ of \mref{famsys} in $X$  {\bf uniformly} (also in $t\in (0,T)$) satisfies
\beqno{contcond} \liminf_{R\to0}\| W\|_{BMO(\Og_R)}=0,\eeq \beqno{DWM}\|W\|_{L^\infty(Q)}+\itQ{Q}{|DW|^{2q}}\le M.\eeq
Then \mref{fullsys0a} has a unique strong weak solution $W$. Moreover, $W\in L^\infty((0,T),W^{1,2p}(\Og))$  if $1-1/p<\nu_*$. In fact, this strong weak solution is a (unique) strong one.\etheo

For $N\le 3$, as corollaries of the above theorem, we can relax the conditions of the them by requiring that solutions to \mref{famsys} are uniformly bounded and continuous. If $N=2$ we merely need that these solutions are uniformly bounded. Moreover, in some cases, the boundedness of their $BMO$ norms will be sufficient to obtain the same results.

In fact, our system \mref{fullsys0a} is inspired by the following   model of two equations ($m=2$) in population biology introduced by Shigesada {\it et al.} in \cite{SKT} has been extensively studied in the last few decades
\beqno{e0a}\left\{\barr{lll} u_t &=& \Delta(d_1u+\ag_{11}u^2+\ag_{12}uv)+k_1u+\bg_{11}u^2+\bg_{12}uv,\\v_t &=& \Delta(d_2v+\ag_{21}uv+\ag_{22}v^2)+k_2v+\bg_{21}uv+\bg_{22}v^2.\earr\right.\eeq   Here, $d_i,\ag_{ij},\bg_{ij}$ and $k_i$ are constants with $d_i>0$. Dirichlet or Neumann boundary conditions were usually assumed for \mref{e0a}. This model (which will be referred to as (SKT) later on) was used to describe the population dynamics of {\em two} species densities  $u,v$ which move and react under the influence of population pressures. Under appropriate assumptions on $d_i,\ag_{ij},\bg_{ij}$ and $k_i$ the global existence of nonnegative solutions with nonnegative data in $W^{2,2}(\Og)$ was established by Yagi in \cite{yag}.

As an extension of these results, we apply our corollaries to a general version of the (SKT) \mref{e0a} which consists of $m$ equations  written compactly as follows
\beqno{skt0a} \left\{\barr{ll}u_t= \Delta (P(u)) +\mg(u)&\mbox{in $\Og\times(0,T_0)$,}\\\mbox{Dirichlet or Neumann boundary conditions}&\mbox{on $\partial\Og\times(0,T_0)$,}\\u=u_0&\mbox{on $\Og$}.\earr\right.\eeq
Here, $P(u)=[P_i(u)]_{i=1}^m$ and $\mg=[\mg_i]_{i=1}^m$ whose components $P_i(u)$'s and $\mg_i$'s
can be appropriate {\em quadratics} in $u\in\RR^m$ (so that the ellipticity condition \mref{ellcond} is satisfied).
This system is the special case of our \mref{fullsys0a} with $\ma =P_u$. In this case, it is obviuous that $\sup_{v\in\RR^m}\frac{|a_v|^2}{\llg(v)}$ is bounded. We will prove that our existence results hold if the norm $\|v\|_{BMO(B_R)}$ is small for $R>0$ is sufficiently small for any strong weak solution $v$. The latter is definitely true for \mref{skt0a} on planar domains as we can control $\|Dv\|_{L^2(\Og)}$.  Thus, the global existence of \mref{skt0a} is completely solved in this situation when $N=2$. Thus, the global existence of \mref{skt0a} is completely solved in this situation.

In \refsec{uni} we prove the uniqueness of strong weak solutions.  For the bounds of fixed points of $\mL$ in $X$, \refsec{Dest} devotes to the estimates of the derivatives of a solution of \mref{fullsys0a} via an use of the Gagliardo-Nirenberg inequality involving BMO norms. In \refsec{proofmain} we will present the proof of \reftheo{main}. Its corollaries and an application to the SKT system \mref{skt0a} will be presented in \refsec{coro} and conclude our paper.

\section{The uniqueness:} \label{uni}\eqnoset

In this section, we prove that a strong weak solution of \mref{fullsys0a} if it exists then it is unique.

Let $W_1, W_2$ be two weak solutions on $Q^t=\Og\times(0,t)$. Set $W=W_1-W_1$. We have 
$$\barrl{W_t=\Div(\ma(W_1)DW_1-\ma(W_2)DW_2)+\Div(\hat{\mb}(W_1)-\hat{\mb}(W_2))+}{3cm}&\hat{\mb}(W_1)DW_1-\hat{\mb}(W_2)DW_2+\mg(W_1)-\mg(W_2).\earr$$
We write
$$\barrl{W_t=\Div(\ma(W_1)DW+[\ma(W_1)-\ma(W_2)]DW_2)+\Div(\hat{\mb}(W_1)-\hat{\mb}(W_2))+}{3cm}&\mb(W_1)DW+[\mb(W_1)-\mb(W_2)]DW_2+\mg(W_1)-\mg(W_2).\earr$$

As $\ma,\hat{\mb},\mb,\mg$ are Lipschitz, we test the above with $W$ to get
$$\barrl{\frac{d}{dt}\iidx{\Og}{|W|^2}+\itQ{Q_t}{|DW|^2}\le}{2cm}& C\itQ{Q_t}{(|W||DW_2||DW|+|W||DW|+|W|^2|DW_2|+|W|^2)}.\earr$$
By Young's inequality this implies
$$\frac{d}{dt}\iidx{\Og}{|W|^2}+\itQ{Q_t}{|DW|^2}\le C\itQ{Q_t}{(|W|^2|DW_2|^2+|W|^2+|W|^2|DW_2|)}.$$
H\"older's inequality yields
$$\barrl{\frac{d}{dt}\iidx{\Og}{|W|^2}+\itQ{Q_t}{|DW|^2}\le}{2cm}& C\dspl{\int_0^t}\left(\iidx{\Og}{|W|^{2q'}}\right)^\frac{1}{q'}\left(\iidx{\Og}{|DW_2|^{2q}}\right)^\frac{1}{q}dt+ C\itQ{Q_t}{|W|^2}.\earr$$

Because  $q>N/2$. So that $2q'<2N/(N-2)$, it is well known that for any $\eg>0$ we have $C(\eg)$ such that
$$\left(\iidx{\Og}{|W|^{2q'}}\right)^\frac{1}{q'}\le \eg\iidx{\Og}{|DW|^2}+C(\eg)\iidx{\Og}{|W|^2}.$$
Combining with the facts that $\|DW_2\|_{L^{2q}(\Og)}\le M$, 
we derive for sufficiently small $\eg$
$$\frac{d}{dt}\iidx{\Og}{|W|^2}+\itQ{Q_t}{|DW|^2}\le C\itQ{Q_t}{|W|^2}.$$

This is a Gr\"onwall inequality which implies $W\equiv0$. Thus, $W_1\equiv W_2$.

\section{Estimates for derivatives} \label{Dest}\eqnoset
In this section, we address the (unniform) bound of fixed points of the map $\mathbb{L}$ by an use of the Gagliardo-Nirenberg inequality involving BMO norms.

We first recall the following simple parabolic version of the usual  Sobolev inequality 
\blemm{parasobo} Let $r=2/N$ if $N>2$ and $r\in(0,1)$ if $N\le 2$. If $g,G$ are sufficiently smooth then 
$$\itQ{\Og\times I}{|g|^{2r}|G|^2}\le C\sup_I\left(\iidx{\Og}{|g|^2}\right)^r\left(\itQ{\Og\times I}{(|DG|^2+|G|^2)}\right).
$$

If $G=0$ on $\partial\Og$ then we can drop the integrand $|G|^2$ on the right hand side. In particular, if $g=G$ we have
$$\itQ{\Og\times I}{|g|^{2(1+r)}}\le C\sup_I\left(\iidx{\Og}{|g|^2}\right)^r\left(\itQ{\Og\times I}{(|Dg|^2+|g|^2)}\right).
$$

\elemm

In the sequel, for $h\ne0$ we will use the difference operator of a function $u=[u_i]_{i=1}^m$
$$\dg_h u(x)=h^{-1}[u(x+he_1)-u(x),\cdots,u(x+he_m)-u(x)].$$

\blemm{Dup} Let $v$ be a strong-weak solution of \mref{famsys} on $\Og\times(0,T)$. Assume that  $v$ satisfies uniformly for any $x_0,t\in Q=\Og\times(\dg_0/2,T)$ that \beqno{contconda} \liminf_{R\to0}\| v\|_{BMO(\Og_R)}=0,\eeq  and  \beqno{DvLq}\|Dv\|_{L^{2q+2}(Q)}<\infty\eeq for  some $q>N/2$. 

Then for some $R_0$ small (depending on the continuity of $v$ in \mref{contconda}) and {\em any} $p$ such that $2p>1$ and $\nu_*>1-1/p$ ({\em that is the spectral gap condition holds for $p$}) then \beqno{Dupest1}\barrl{\sup_{t\in(\dg_0/2,T)}\iidx{\Og}{|Dv|^{2p}}+\itQ{Q}{|Dv|^{2p-2}|D^2v|^2}\le}{2cm}&  C(R_0^{-2},p)\itQ{Q}{|Dv|^{2q}}+\iidx{\Og\times\{\dg_0/2\}}{|Dv|^{2p}}.\earr\eeq

In addition (one should note the exponent  $q$ on the right hand side), $$\|Dv\|_{L^{2p\cg}(Q)} \le C\left(R_0^{-2},\itQ{Q}{|Dv|^{2q}}, \iidx{\Og\times\{\dg_0/2\}}{|Dv|^{2p}}\right).$$

\elemm

\bproof We need only to look at the case $\tau=1$ as the argument is similar for $\tau\in[0,1]$.
Apply $\dg_h$ to  the equation of $v$ to see that $v$ weakly solves
$$(\dg_hv)_t=\Div(\ma D(\dg_h v)+\ma_v \myprod{\dg_h v,Dv})+\Div(\hat{\mb}_v\dg_hv))+\mb_v\myprod{\dg_hv,Dv}+\mb D(\dg_hv)+\mg_v  \dg_hv.$$

Since the parameters $\ma,\hat{\mb}, \mb$ and $\mg$ of the equation are bounded (see \refrem{Bgrowth}), 
for any $0<s<t<2R_0$ we test this system with $|\dg_hv|^{2p-2}\dg_hv\fg^2\eta$, with $p\ge1$ and   $\fg,\eta$ being positive $C^1$ cutoff functions for the concentric balls $B_s, B_t$ and the time interval $I$, and use Young's inequality for the term $|\ma_v||\dg_hv|^{2p-1}|Dv||D(\dg_hv)|$ and the spectral gap condition (with $X=\dg_h v$) to get  (we refer to \cite{dleJMAA} for details). We see that  for some constant $C$ and $Q=\Og_t\times (\dg_0/2,T)$
$$\barrl{\sup_{t\in(\dg_0/2,T)}\iidx{\Og}{|\dg_hv|^{2p}\fg^2}+\itQ{Q}{ |\dg_hv|^{2p-2}|D(\dg_hv)|^2\fg^2}\le C\itQ{Q}{|\dg_hv|^{2p}|Dv|^2\fg^2}}{2cm}&+C\itQ{Q}{(|\dg_hv|^{2p}+|Dv|^{2p})(\fg^2+|\dg_h\fg|^2+|D\fg|^2)}+\iidx{\Og\times\{\dg_0/2\}}{|\dg_hv|^{2p}\fg^2}.\earr$$
Another use of Young's inequality for the first term on the right hand side yields
$$\barrl{\sup_{t\in(\dg_0/2,T)}\iidx{\Og}{|\dg_hv|^{2p}\fg^2}+\itQ{Q}{ |\dg_hv|^{2p-2}|D(\dg_hv)|^2\fg^2}\le C\itQ{Q}{|\dg_hv|^{2p+2}\fg^2}+ }{2cm}&C\itQ{Q}{|Dv|^{2p+2}\fg^2}+C\itQ{Q}{(|\dg_hv|^{2p}+|Dv|^{2p})(\fg^2+|\dg_h\fg|^2+|D\fg|^2)}+\\&\iidx{\Og\times\{\dg_0/2\}}{|\dg_hv|^{2p}\fg^2}.\earr$$

For any $p\ge1 $ such that (this is true if $p=q$ because of \mref{DvLq})
\beqno{keyassDv}\itQ{Q}{|Dv|^{2p+2}}<\infty\eeq we let $h\to0$ and see that
\beqno{keyassDv1}\barrl{\sup_{t\in(\dg_0/2,T)}\iidx{\Og}{|Dv|^{2p}\fg^2}+\itQ{Q}{|Dv|^{2p-2}|D^2v|^2\fg^2}\le}{.5cm}& C\itQ{Q}{|Dv|^{2p+2}\fg^2}+ C\itQ{Q}{|Dv|^{2p}(\fg^2+|D\fg|^2)}+\iidx{\Og\times\{\dg_0/2\}}{|Dv|^{2p}\fg^2}.\earr\eeq

The key issue here is that we will have to handle  the integral of $|Dv|^{2p+2}$.

For any $t>0$ we write $\Og_t=\Og\cap B_t$, $Q_t=\Og_t\times(\dg_0/2,T)$ and $$\ccA_p(t)=\sup_{(\dg_0/2,T)}\iidx{\Og_t}{|Dv|^{2p}},\;\ccB_p(t)=\itQ{Q_t}{|Dv|^{2p+2}},f(t)=\iidx{\Og_t\times\{\dg_0/2\}}{|Dv|^{2p}},$$
$$\ccH_p(t)=\itQ{Q_t}{|Dv|^{2p-2}|D^2v|^2},\; \ccG_p(t)=\itQ{Q_t}{|Dv|^{2p}}.$$

By the local Gagliardo-Nirenberg inequality (\cite[Lemma 2.4]{dleANS} with $u=U$, $\Fg\equiv1$ and $\psi$ is a cutoff function for $B_s, B_t$) for any $\eg>0$ and some constant $C=C(\eg,N)$ we have that (in this paper we will refer to this as the Gagliardo-Nirenberg BMO inequality) $$\barr{lll}\iidx{\Og_s}{|Dv|^{2p+2}}&\le& \eg \iidx{\Og_t}{|Dv|^{2p+2}}+
\\&&
C\|v\|_{BMO(\Og_t)}^2\iidx{\Og_t}{|Dv|^{2p-2}|D^2v|^2}+C\frac{\|v\|^2_{BMO(\Og_t)}}{(t-s)^2}\iidx{\Og_t}{|Dv|^{2p}}.\earr$$

Because of \mref{keyassDv} and \mref{keyassDv1}, the terms in this inequality are all finite and it holds for a.e. $t$. Using this inequality in \mref{keyassDv1}, where $\fg$ is a cut-off function for $B_s,B_t$ and by the assumption \mref{contconda} on the uniform continuity of $v$, $\|v\|_{BMO(\Og_t)}$ can be very small. For $0<s<t<R_0$ with sufficiently small $R_0$ depending on the uniform continuity of $v$ we obtain the following recursive system of inequalities
$$\ccA_p(s)+\ccH_p(s)\le C\ccB_p(t)+ \frac{C}{(t-s)^2}\ccG_p(t)+f(t),$$
$$\ccB_p(s)\le \eg (\ccH_p(t)+\ccB_p(t))+\frac{C}{(t-s)^2}\ccG_p(t).$$
If $\eg$ (or $R_0$) is small then we can iterate this to absorb the terms $\ccH_p,\ccB_p$ on the right hand side to the left hand side to get (see \cite[inequality (3.27), proof of Proposition 3.1]{dleANS})
\beqno{Duiter}\ccA_p(R_0)+\ccH_p(R_0)\le \frac{C}{R_0^2}\ccG_p(2R_0)+f(R_0).\eeq

This gives the estimate of $\ccA_p, \ccH_p$. That is
$$\barrl{\sup_{(\dg_0/2,T)}\iidx{\Og_{R_0}}{|Dv|^{2p}}+\itQ{Q_{R_0}}{|Dv|^{2p-2}|D^2v|^2}\le}{2cm}& CR_0^{-2}\itQ{Q_{2R_0}}{|Dv|^{2p}}+\iidx{\Og\times\{\dg_0/2\}}{|Dv|^{2p}}.\earr$$

This also holds when $B_{R_0}$ intersects the boundary $\partial\Og$ see \cite[Remarks 3.3.5 and 3.3.6]{dlebook}. Fixing such $R_0$ and covering $\Og$ with balls of radius $R_0$ and summing the above inequality over this partition, we derive a global estimate
\beqno{DuLp}\barrl{\sup_{(\dg_0/2,T)}\iidx{\Og}{|Dv|^{2p}}+\itQ{\Og\times(\dg_0,T)}{|Dv|^{2p-2}|D^2v|^2}\le}{2cm} &C(R_0^{-2})\itQ{\Og\times(\dg_0,T)}{|Dv|^{2p}}+\iidx{\Og\times\{\dg_0/2\}}{|Dv|^{2p}}.\earr\eeq

Dropping the variable $R_0$ and $W_0$ (it is the value of a classical solution), we define $\ccA_p,\ccH_p,\ccB_p$ in the same ways with $\Og_{R_0}=\Og$. By the  parabolic Sobolev inequality in \reflemm{parasobo} with $g=G=Dv$, \mref{Duiter} also shows that its right hand side is self-improved.  That is if $\ccG_p$ is finite for some $p\ge1$ then so are $\ccA_p, \ccH_p$. We have
$$\|Dv\|_{L^{2p\cg}(Q)} \le C\left(R_0^{-2},\itQ{Q}{|Dv|^{2p}}, \iidx{\Og\times\{\dg_0/2\}}{|Dv|^{2p}}\right),$$ where $\cg= 1+2/N$ if $N\ge3$ and $\cg=1+r$ for any $r\in(0,1)$ if $N=2$. 
Thus, \mref{keyassDv} holds again with $2p+2$ is now $2p\cg$. This also establishes our last assertion of the lemma (keep in mind that the values of $v$ in $(0,\dg_0)$ are those of the strong solution).

This argument can be repeated with $p$ being replaced by $\tau p$ for $\tau>1$ as long as $2\tau p+2\le \cg 2p$ (see \mref{keyassDv}). 
Define $p_{n+1}=\tau_np_n$ with $\tau_n=\cg-1/p_n$. We see that $2p_{n+1}+2=2p_n\cg$. Moreover, $\tau_n> 1\Leftrightarrow p_n>N/2\Leftrightarrow p_n\uparrow$. Thus, if we take $p_1$ to be the number $q>\frac{N}{2}$ in \mref{DvLq} of this lemma then such sequences $\{p_n\},\{\tau_n\}$ exist and the iterate the argument as long as $\nu_*>1-1/p_n$.

Along the sequence $\{p_n\}$ we have (if the spectral gap condition holds for $p_n$)
$$\barrl{\sup_{(\dg_0/2,T)}\iidx{\Og}{|Dv|^{2p_n}}+\itQ{Q}{|Dv|^{2p_n-2}|D^2v|^2}\le}{4cm}& C\left(n,R_0^{-2},\itQ{Q}{|Dv|^{2p_{1}}},\iidx{\Og\times\{\dg_0/2\}}{|Dv|^{2p}}\right).\earr$$

 As $p_n\to\infty$ since $\tau_n>1$ and $p_1=q$, we then have \mref{Dupest1}. 
\eproof

\brem{Bgrowth} We make use of the polynomial growth of $ \mb,\mg$ in deriving \mref{keyassDv1}. 
\erem

\section{The proof of the main result}\label{proofmain}{\eqnoset}

In this section, we will prove the existence of fixed points of $\mL$ in an appropriate defined space $X$ via the Leray-Schauder theory.

First of all, we address the compactness of the operator $\mL$ that leads to the definition later.

\subsection{A compactness lemma}

The main issues is to show $\mL(K)$ of a bounded set $K$ of $X$ is precompact. The continuity of $\mL$ follows in a standard way. 
Similar to \cite[Lemma 3.3]{letrans}, we will establish the following compactness result which will serve the purpose.

\blemm{complem} Let $\ccF$ be a collection of function $w$ such that $D^2w$ exists and satisfies
\beqno{vccF}\left|\itQ{Q}{Dw\psi_t}\right|\le C\itQ{Q}{(|D^2w||D\psi|+|Dw||\psi|)} \quad \forall\psi\in C^1(Q)\eeq
and $\psi(\cdot,0)\equiv\psi(\cdot,T)\equiv0$. Suppose that for some constants $M$, $p\ge1$ and  for all $w\in\ccF$ \beqno{Mwbound} \sup_{t\in (0,T)} \iidx{\Og}{|Dw|^{2p}}+\itQ{Q}{|Dw|^{2p-2}|D^2w|^2} <M,\eeq
\beqno{Mwbound1}\itQ{Q_{s,t}}{|D^2w|^2} =O(t-s),\eeq
where $Q_{s,t}=\Og\times(s,t)$.
Then 

\bdes\item[i)]$\{Dw\,:\,w\in\ccF\}$ is compactly embedded in $L^{2p}((0,T),L^{2p}(\Og))$.

\item[ii)] Let $\cg_0=1+2/N$. If $ \cg\in (1,\cg_0)$ then $\{Dw\,:\,w\in \ccF\}$  is compactly embedded in $L^{2p\cg}((0,T),L^{2p\cg}(\Og))$.
\edes

\elemm

\bproof  For each $w\in\ccF$ denote $v=|Dw|^{p-1}Dw$. We first show that for $l>(N+2)/2$ \beqno{wcont}\|v(\cdot,t+h)-v(\cdot,t)\|_{W^{-l,2}(\Og)}\le C[M\eg_0 +C(\eg_0)O(h)]\quad \forall \eg_0>0.\eeq

For any $Q_{s,t}=\Og\times(s,t)$ let $h=t-s$ and $\psi(x,t)=|Dw|^{p-2}Dw\eta(t)\fg(x)$ in \mref{vccF} with $\fg\in C^1(\Og)$ and $\eta\equiv1$ in $(s,r)$ and $\eta\equiv0$ outside $(s-\eg,r+\eg)$. From \mref{vccF}, we get
$$\left|\itQ{Q_{s,t}}{v\psi_t}\right|\le C\itQ{Q_{s,t}}{|\eta(|D^2w|^2|Dw|^{p-2}\fg +|D^2w||Dw|^{p-1}D\fg|+|Dw|^p|\fg|)|}.$$

We estimate the right hand side. By Young's inequality, for every $\eg_0>0$ we have $$\itQ{Q_{s,t}}{\eta|D^2w|^2|Dw|^{p-2}\fg}\le \itQ{Q_{s,t}}{|D^2w|^2(\eg_0|Dw|^{2p-2}+C(\eg_0))}\|\fg\|_{C^1(\Og)},$$
$$\itQ{Q_{s,t}}{\eta|D^2w||Dw|^{p-1}D\fg}\le |h|^\frac12\left(\itQ{Q_{s,t}}{|D^2w|^2|Dw|^{2p-2}}\right)^\frac12\|\fg\|_{C^1(\Og)}$$
$$\itQ{Q_{s,t}}{\eta|Dw|^{p}\fg}\le |h|^\frac12\left(\itQ{Q_{s,t}}{|Dw|^{2p}}\right)^\frac12\|\fg\|_{C^1(\Og)}.$$

Because $l>(N+2)/2$, $\|\fg\|_{C^1(\Og)}\le C\|\fg\|_{W^{l,2}(\Og)}$, by the assumptions \mref{Mwbound} and \mref{Mwbound1} we see that the above inequalities together imply
\beqno{keycomplemm}\left|\itQ{Q_{s,t}}{v\psi_t}\right|\le C[M\eg_0 +C(\eg_0)O(h)]\|\fg\|_{W^{l,2}(\Og)} \quad \forall\fg\in W^{l,2}(\Og).\eeq

Letting $\eg\to0$, we also have  for all $\fg\in W^{l,2}(\Og)$ (compare with \cite[Lemma 3.2]{dletrans})
$$\left|\iidx{\Og}{[v(\cdot,s)-v(\cdot,t)]\fg}\right|=\left|\itQ{Q_{s,t}}{v_t\psi}\right|\le C[M\eg_0 +C(\eg_0)O(h)]\|\fg\|_{W^{l,2}(\Og)}.$$

Hence \mref{wcont} follows.  We now interpolate $L^2(\Og)\cap W^{1,2}(\Og)$ between $W^{1,2}(\Og)$ and $W^{-l,2}(\Og)$ to get for any $\mu>0$ and $v\in L^2(\Og)\cap W^{1,2}(\Og)$ that
$$\|v(\cdot,t+h)-v(\cdot,t)\|_{L^2(\Og)} \le \mu \|v\|_{W^{1,2}(\Og)}^2+C(\mu)\|v(\cdot,t+h)-v(\cdot,t)\|_{W^{-l,2}(\Og)}^2$$

Since $\dspl{\int_0^T}\|v\|_{W^{1,2}(\Og)}^2dt\le C(M)$ if $v\in \ccF$, for any given $\eg>0$, we can choose $\mu$ small first and then $\eg_0,h$ small such that by \mref{wcont}
$$\int_{0}^{T-h}\|v(\cdot,t+h)-v(\cdot,t)\|_{L^2(\Og)}^2dt<\eg\quad \forall v\in\ccF.$$
Thus, we just prove that
$$\int_{0}^{T-h}\|v(\cdot,t+h)-v(\cdot,t)\|_{L^2(\Og)}^2dt\le O(h).$$

For any $t_1,t_2\in (0,T)$ and $v\in \ccF$, the fact that the collection of $\int_{t_1}^{t_2}vdt$ is a pre-compact set of $L^2(\Og)$ is clear because the set $$\{ \frac{1}{|t_2-t_1|}\int_{t_1}^{t_2}vdt\,:\,  \int_{t_1}^{t_2}\iidx{\Og}{\Fg(|v|^2+|Dv|^2)}dt <M)\}$$ belongs to the closure of the convex hull in $L^2(\Og)$, a bounded set in $W^{1,2}(\Og)$ as $\Fg\ge \llg_0>0$, which is compact in $L^2(\Og)$.

By a result of Simon \cite{JS} as in \cite[Lemma 3.3]{letrans} (setting $\ccB=L^2(\Og)$), we proved the compactness of $\{v:w\in\ccF\}$ in $L^2(0,T,L^2(\Og))$. 

By the definition of $v$, the set $\{{|Dw|^{p-1}Dw\,:\, w\in \ccF}\}$ is compact in $L^2(0,T,L^2(\Og))$. It is well known that if a sequence $\{f^{2p}_n\}$ converges in $L^1(\Og)$ then $\{f_n\}$ converges in $L^{2p}(\Og)$ (by dominating convergence theorem). This implies that $\{Dw\,:\,w\in\ccF\}$ is compact in $L^{2p}(0,T,L^{2p}(\Og))$. This gives i).

Finally, if $\cg\in(1,\cg_0)$ then  there are  $\ag,\bg \in (0,1)$ such that $2p\cg=2p\ag+2p\cg_0\bg$. 
Let $\{Dw_n\}$ with $w_n\in\ccF$ be a bounded set satisfying \mref{Mwbound} and \mref{Mwbound1}. By i) and relabeling, we can assume $\{Dw_n\}$ is convergent sequence in $L^{2p}(0,T,L^{2p}(\Og))$. We write $W_{n,m}= w_n-w_m$. By H\"older's inequality
$$\barr{lll}\itQ{Q}{|DW_{n,m}|^{2p\cg}}&\le&\left(\itQ{Q}{|DW_{n,m}|^{2p}}\right)^\frac{1}{\ag}\left(\itQ{Q}{(|DW_{n,m}|)^{2p\cg_0}}\right)^\frac{1}{\bg}\\&\le&\left(\itQ{Q}{|DW_{n,m}|^{2p}}\right)^\frac{1}{\ag}\left(\itQ{Q}{(|Dw_n|+|Dw_m|)^{2p\cg_0}}\right)^\frac{1}{\bg}.\earr$$
The first factor on the right goes to zero because $\{DW_{n,m}\}$ converges to 0 in $L^{2p}(0,T,L^{2p}(\Og))$.  By \mref{Mwbound} and the parabolic Sobolev inequality in \reflemm{parasobo} with $|g|=|G|=|Dw_n|^p$, $\{Dw_n\}$ and $\{Dw_m\}$ are bounded uniformly in $L^{2p\cg_0}(0,T,L^{2p\cg_0}(\Og))$. The second factor is bounded. Thus $\ccF$ is compact in $L^{2p\cg}(0,T,L^{2p\cg}(\Og))$. This gives ii). The proof is complete. \eproof

\subsection{The space $X$ and $\mL:X\to X$ is compact}

Next, we will define the space $X$ such that $\mL\,:\, X\to X$ is a compact. By the theory of \cite{Am2} the linear parabolic system defining $\mL(\Fg)$ has sufficient smooth coefficients so that it has  a classical solution. However, in what below we will need some uniform estimates to show that $\mL :X\to X$ is a compact map. That is $\mL(K)$ is compact in $X$ if $K\subset X$ is a bounded set.

Apply $\dg_h$ to  the equation of $W$ to see that 
\beqno{difhW}\barrl{(\dg_hW)_t=\Div(\ma(\Fg) D(\dg_h W)+\ma_\Fg \myprod{\dg_h \Fg,DW})+\Div(\hat{\mb}_\Fg\dg_h\Fg)+}{4cm}&\mb_\Fg\dg_h\Fg DW+\mb D\dg_h W+\mg_\Fg \dg_h\Fg.\earr\eeq

Assume that the condition GS holds for $q$. That is $\nu_*>1-\frac1q$. We test this system with $|\dg_h W|^{2q-2}\dg_h W$  and let $h\to 0$ to have
$$\barrl{\iidx{\Og}{|DW|^{2q}}+\itQ{Q}{|DW|^{2q-2}|D^2W|^2}\le }{2cm}&
C\itQ{Q}{(|D\Fg||DW||DW|^{2q-2}|D^2W|+|D\Fg||DW|^{2q-2}|D^2W|)}+\\&
 +C\itQ{Q}{(|D\Fg||DW|^{2q-1}+|DW|^{2q-1}|D^2W|+|D\Fg||DW|^{2q-1})}.\earr$$

Applying Young's inequality
\beqno{defX2}\barrl{\iidx{\Og}{|DW|^{2q}}+\itQ{Q}{|DW|^{2q-2}|D^2W|^2}\le }{2cm}&
C\itQ{Q}{(|D\Fg|^2|DW|^{2q}+|D\Fg|^{2q}+|DW|^{2q})}.
\earr
\eeq

Being inspired by the above calculations, \reflemm{Dup}, the compactness result and \reflemm{Lcompact} below, we introduce the space $X$ here such the right hand side of \mref{defX2} is finite so that $W\in X$ (by letting $h\to0$). 

Assume that for some $q_0>N/2$ with $\nu_*>1-\frac{1}{q_0}$.  Let $\cg\in(1,1+2/N)$ be such that  $q_0(\cg-1)>1$.

Since $q_0(\cg-1)>1$, we see that $\cg-1/q_0>1$ so that we can fix a number $\ag\in(0,1)$ such that  $1<\ag\cg\le \cg-1/q_0$. That is $\ag\cg>1$ and $\frac{2}{1-\ag}\le 2q_0\cg$. 
We summarize the choices of $q_0,\cg,\ag$ below
\beqno{GSi}\nu_*>1-\frac{1}{q_0},\; \ag\cg>1,\;\frac{2}{1-\ag}\le 2q_0\cg.\eeq

For some $\eg_0>0$ define $X=\{\Fg\,:\,D\Fg\in L^{2q_0\cg}(Q), \Fg\in C^{0,\eg_0}(Q)\}$ with norm
$$\|\Fg\|_X=\|\Fg\|_{C^{0,\eg_0}(Q)}+\|D\Fg\|_{L^{2q_0\cg}(Q)}.$$ 
Of course, we can choose $\cg\sim 1+2/N$ and $\ag\sim 1/\cg$ so that $q_0\sim 1/(\cg-1)\sim N/2$.

To begin, we have the following simple result.

\blemm{lsu} If $\Fg\in X$ then $W=\mL(\Fg)\in W^{2,1}(Q)$ and we have the following estimate $$\|W\|^{(2)}_{2,Q}\le C(\|\Fg\|^{(1)}_{2,Q}+\|W_0\|^{(1)}_{2,Q}).$$ Here, $\|W\|^{(l)}_{2,Q}=\sum_{2i+j=l}\|D_t^i D_x^jW\|_{L^2(Q)}$ and the constant $C$ depends on $\|\Fg\|_{L^\infty(Q)}$.
\elemm
\bproof
We can rewrite the system \mref{famsys} as
${\cal L}(W)=f$ with
$$ {\cal L}(W) = W_t-\ma D^2W-[\ma_\Fg(\Fg)D\Fg+\mb(\Fg)]DW  \mbox{ and } f=\Div(\hat{\mb}(\Fg))+\mg(\Fg).$$

The assertion is a simple consequence of \cite[Theorem 9.1]{LSU} which can be extended to linear systems with smooth coefficients by a similar study of fundamental solutions for systems. We check the conditions of \cite[Theorem 9.1]{LSU}. By the definition of $X$, $\ma$ is continuous on $Q$ and bounded. Moreover, because $q_0>N/2$ and $\cg\in(1,1+2/N)$, we can choose $\cg$ such that $r=2q_0\cg>N+2$ and $\ma(\Fg)D\Fg+\mb(\Fg)\in L^{loc}_{r,Q}$. We also have $\mg(\Fg)\in L^{loc}_{s,Q}$ for any $s\in(1,\infty)$. Thus, the assumption \cite[(7.1) of Theorem 9.1]{LSU} is verified.

Next, since the initial value $W_0$ is the value of $W_*(\dg_0/2)$ of the classical solution $W_*$, the compatibility condition between the initial and boundary data ($\fg=W_0$ and $\Fg=0$ in \cite{LSU}) of \cite[(9.2)]{LSU} holds.

Thus, with $q=2$, we have from \cite[Theorem 9.1]{LSU} that there is a constant $C$ depending on $\|\Fg\|_{L^\infty(Q)}$ such that for $\|f\|_{2,Q}=\|f\|_{L^2(Q)}$
$$\|W\|^{(2)}_{2,Q}\le C(\|f\|_{2,Q}+\|W_0\|^{(1)}_{2,Q}),$$
From the definition of $f$, this completes the proof. \eproof

\brem{SLSU}Moreover, if $\Fg\in X$ then $f\in L^{2q_0\cg}(Q)$ with $2q_0\cg>(N+2)/2$. By \cite[Theorem 2.1 of Chapter VII]{LSU} and the system defines $W$ has continuous bounded coefficients, we see that $\|W\|_{L^\infty(Q)}$ is bounded. 
\erem

We now apply the above argument to show that
\blemm{Lcompact} 
Assume \mref{GSi}. If  $\eg_0>0$ is suffiently small then $\mL:\,X\to X$ is compact. Moreover, $\mL(\Fg)\in W^{1,2q_0}(\Og)$.
\elemm
\bproof Let $W=\mL(\Fg)$. For some  fixed $R_0>0$ and $Q_{R_0}=\Og_{R_0}\times(0,T)$
By the Gagliardo-Nirenberg BMO inequality, we have $$\itQ{Q_{R_0}}{|DW|^4}\le C\|W\|^2_{L_\infty(Q)}\itQ{Q_{2R_0}}{|D^2W|^2}+C(R_0^{-2})\|W\|^2_{L^\infty(Q)}\itQ{Q_{2R_0}}{|DW|^2}.$$

From \refrem{SLSU} if $\Fg\in X$ then $W$ is bounded in terms of $\|\Fg\|_X$. So that by summing over a finite covering of $\Og$ by  balls of radius $R_0$, we obtain
$$\itQ{Q}{|DW|^4}\le C(R_0,\|\Fg\|_X)\itQ{Q}{|D^2W|^2}+C\itQ{Q}{|DW|^2}.$$

By testing the system with $W$, we easily see that $\|DW\|_{L^2(Q)}$ is bounded in terms of $\|\Fg\|_X$. Also, by \reflemm{lsu}  $\|D^2W\|_{L^2(Q)}$ is bounded in terms of $\|\Fg\|_X$
We conclude that $\|DW\|_{L^4(Q)}$ is bounded in terms of $\|\Fg\|_X$. That is, for some any $\cg\in(1,1+2/N)$
$$\itQ{Q}{|DW|^{2\cg}}<\infty.$$

Let $q_1=1$.
For some $\ag\in(0,1)$ such that $\ag\cg>1$ we define $q_{i+1}=q_i\ag\cg$. Let $q=q_{i+1}$ in \mref{defX2}. If the condition GS holds for $s=2q_{i+1}-2$ (that is $\nu_*>1-1/q_{i+1}$) then 
$$\iidx{\Og}{|DW|^{2q_{i+1}}}+\itQ{Q}{|DW|^{2q_{i+1}-2}|D^2W|^2}\le \itQ{Q}{(|D\Fg|^{2}|DW|^{2q_{i+1}}+1)}.$$
If $DW\in L^{2q_i\cg}(Q)$,  $i\ge1$,
we apply Young's inequality to the term $|D\Fg|^{2}|DW|^{2q_{i+1}}$ on the right hand side to obtain
\beqno{DWq}\iidx{\Og}{|DW|^{2q_{i+1}}}+\itQ{Q}{|DW|^{2q_{i+1}-2}|D^2W|^2}\le \itQ{Q}{(|D\Fg|^{\frac{2}{1-\ag}}+|DW|^{\frac{2q_{i+1}}{\ag}}+1)}.\eeq

Because $2q_0\cg\ge\frac{2}{1-\ag}$, $D\Fg\in L^{\frac{2}{1-\ag}}(Q)$. Also, since $DW\in L^{2q_i\cg}(Q)$ and $2q_i\cg=\frac{2q_{i+1}}{\ag}$,  so that the above quantities are all finite. Again, by parabolic Sobolev inequality, we have $DW\in L^{2q_{i+1}\cg}(Q)$. We can repeat the argument as it is true for $q_1$ to see that $DW\in L^{2q_i\cg}(Q)$ for all $i\ge1$ (as long as the spectral gap condition holds for $q_i$). Note that $q_i\to\infty$, because $\ag\cg>1$.

Under the assumption that the condition GS that $\nu_*>1-1/q_i$, which holds for any exponent $q_i$ such that $q_i\le q_0$ because of \mref{GSi}, we have $DW\in L^{2q_i\cg}(Q)$ by \mref{DWq} and the parabolic Sobolev inequality. Thus, $\mL(\Fg)\in X$ (choose $i$ such that $q_i=q_0$).

Now, if $\Fg$ belongs to a bounded set $K$ of $X$ then \mref{DWq} shows that the set  $\mL(K)$ satisfies the conditions of \reflemm{complem}. Therefore, the set $\{DW:W\in \mL(K)\}$ is compact in $L^{2q_i\cg}(0,T,L^{2q_i\cg}(\Og))$. Therefore, $\mL\,:\, X\to L^{2q_0\cg}(0,T,L^{2q_0\cg}(\Og))$ is compact. Note that, by testing with $\psi(\cdot,0)\equiv\psi(\cdot,T)\equiv0$, we have 
$$\left|\itQ{Q}{DW\psi_t}\right|\le C\itQ{Q}{(|D^2W||D\psi|+|DW||D\Fg||D\psi|+(|D\Fg|+1)|DW||\psi|)}.$$

There are some extra terms involving $\Fg\in X$ but, by applying the H\"older or Young inequalities and if $\|\Fg\|_X$ is bounded, we can see that we still get \mref{keycomplemm} in the proof of \reflemm{complem} which is obtained from the assumption \mref{vccF}, namely
$$\left|\itQ{Q}{DW\psi_t}\right|\le C\itQ{Q}{(|D^2W||D\psi|+|DW||\psi|)} \quad \forall\psi\in C^1(Q).$$ 
Thus, \reflemm{complem} is applicable here.

Moreover, we see that $D^2W\in L^2(Q)$  by letting $h\to0$ in \mref{difhW}. If $\Fg$ is in a bounded set of $X$ then $\|D^2W\|_{L^2(Q)}$ is uniform bounded. From the above argument, if $q_0>N/2$ then for some $q\in(N/2,q_0]$ we also see that $\|DW\|_{L^{2q}(\Og)}$ is uniformly bounded. So,  $W$ is H\"older continuous in $x$, as $W\in W^{1,2q}(\Og)$ with $2q>N$.  Since $\|\Fg\|_{L^\infty}(Q)$, $\|D^2W\|_{L^2(Q)}$ and $\|D\Fg\|_{L^{2\cg}(Q)}$ and $\|DW\|_{L^{2\cg}(Q)}$ are uniform bounded, by \mref{difhW} we can solve for $W_t$ and see that $W_t\in L^{2r}(Q)$ with uniform bounded norm for some $r\in(0,1)$. We conclude that $W$ is H\"older continuous in $x,t$ (see \cite[Lemma 4]{NSver}). Thus, $\mL$ is also compact in the space $C^{0,\eg_0}(Q)$ if $\eg_0$ is sufficiently small.

Hence, if $q_0>N/2$ and $\eg_0$ small then $\mL\,:\, X\to X$ is compact and $\mL(\Fg)=W\in W^{1,q_0}(\Og)$. The proof is complete. \eproof

\brem{qrem} If  $W\in X$ then $DW\in L^{2q\cg}(Q)$ for $q\le q_0$. If $q\in(N/2,q_0)$ we also have $DW\in L^{2q+2}(Q)$. Indeed, we have $2q+2<2q\cg$ because this is equivalent to $q>N/2$ for some $\cg\in(1,1+2/N)$ (or $\cg\in(1,2)$ if $N=2$).

\erem

\brem{higherintegrable} 

The higher integrability of $DW$ for $\Fg\in X$ (in particular $DW\in L^4(Q)$, or the results in \cite{GiaS} which can be extended to the parabolic boundary of $Q$ but the number $q_0$ must be redefined in order that the gap condition to be satisfied) is crucial because we could have not started the iteration argument (starting wih \mref{defX2} to get $q_1=1$) without it. On the other hand, if we assumed $D\Fg\in L^\infty(Q)$ but our argument could not provides $DW\in L^\infty(Q)$.\erem

\subsection{Proof of the main theorem}

We see now that \reftheo{main} is proved if we can establish

\blemm{unifp} Assume that \mref{GSi} holds and $q_0>N/2$.
If there is $q\in(N/2,q_0)$ and $M$ such that any strong weak solution $W$ of \mref{famsys} in $X$  {\bf uniformly} (also in $t$) satisfies
\beqno{contcond} \liminf_{R\to0}\| W\|_{BMO(\Og_R)}=0,\eeq \beqno{DWM}\|W\|_{L^\infty(Q)}+\itQ{Q}{|DW|^{2q}}\le M.\eeq
Then \mref{fullsys0a} has a strong weak solution in $X$.
\elemm
{\bf Proof of \reftheo{main} or \reflemm{unifp}:}
Thanks to \reflemm{Lcompact},  $\mL:X\to X$ is a compact map. Let $W$ be a strong weak solution of \mref{famsys}  in $X$.  From the proof of \reflemm{Lcompact} we see that $W\in  W^{1,2q}(\Og)$  and  $W\in  W^{1,2q+2}(Q)$ for any $q\in (N/2,q_0]$ because $2q+2\le 2q\cg$ (see \refrem{qrem}). Hence $W$ is bounded and H\"older continuous because $W\in W^{1,2q}(\Og)$ and $2q>N$ (but this continuity may not be uniform among such $W$ so that we can apply \reflemm{Dup} but the obtained estimates are not {\em uniform} because the number $R_0$ is not fixed).

By the assumption of the Theorem there is $q,M$ such that $q>N/2$ and $$\itQ{Q}{|DW|^{2q}}<M$$  and \mref{contcond} holds  uniformly. So that for some uniform small $R_0$, \reflemm{Dup} yields that  
$$\iidx{\Og}{|DW|^{2p}}+\itQ{Q}{|DW|^{2p-2}|D^2W|^2}\le C\left(p, R_0^{-2},\itQ{Q}{|DW|^{2q}}, \iidx{\Og}{|D W_0|^{2p}}\right)$$
for any $p$ such that $\nu_*>1-1/p$. This condition holds for $p=q_0$ because of \mref{GSi}. One should recall that $W_0(x)$ is $W_*(x,\dg_0)$, the value of the strong solution. Choose $p$ such that $p=q_0$ then $\|DW\|_{L^{2q_0\cg}(Q)}$ is uniformly bounded. The bound for $\|W\|_{C^{0,\eg_0}(Q)}$ is obvious from \reflemm{Lcompact}. Therefore,
we have the uniform bound for such $W$ in $X$. The existence of a strong weak solution then follows. 

The above argument holds for $2p>N$ and yields
$$ \sup_{(0,T_0)}\iidx{\Og}{|DW|^{2p}}<\infty$$ so that our strong weak solution (which is unique) coincides with the strong solution. \eproof
 
\brem{scalareqn} The above argument applies to scalar equation for any dimension $N$ (or systems if $N=2$), the spectral gap condition is not needed. If a fixed point $u$ is H\"older continuous then $Du\in L^{2p}(Q)$ for all $p$.

\erem
\brem{condscont} The uniformity in the condition \mref{contcond} is essential and \mref{DWM} does not implies it even if $q>N/2$. Without this uniform continuity,  we can not obtain \reflemm{Dup} to get uniform bound for fixed points of $\mL$ in $X$ because the number $R_0$ in \mref{Dupest1} of  \reflemm{Dup} is not uniform.
\erem

\section{The corollaries}\label{coro}\eqnoset
When $N\le3$, the condition on $\|DW\|_{L^{2q}(Q)}$ in \mref{DWM} is almost obvious and we have

\bcoro{Nlessequal3} If $N\le 3$ and $\nu_*>1-2/N$. Suppose that any solution of \mref{fullsys0a} is uniformly bounded and uniformly  continuous then there is a unique strong weak solution. \ecoro 

\bproof Note that  \mref{ellcond10} holds for $s>-1$ and $\nu>s/(s+2)$ so with $s=2p-2$ then if $p>1/2$ and $\nu>1-1/p$ then we can still apply the argument of \reflemm{Dup} to obtain a uniform bound for the fixed points of $\mL$ in $X$.

Thus, under the assumption that all weak solutions of \mref{fullsys0a} are uniformly bounded and continuous we see that \mref{contcond} holds uniformly so that we can start with $p>1/2$ and assume that $\nu>1-1/p$ to prove that (see \reflemm{Dup}): if $v$ is a fixed point of $\mL$ in $X$ then for $\cg=1+2/N$ the following integrability improvement holds $$ Dv\in L^{2p}(Q) \Rightarrow Dv\in L^{2p}(\Og) \mbox{ and } |Dv|^{p-1}D^2v\in L^2(Q)\Rightarrow Dv\in L^{2p\cg}(Q).$$

We see that $2p\cg >N$ if $p>\frac{N^2}{2(N+2)}$. Thus, we need $\frac{N^2}{2(N+2)}>\frac12$, which is equivalent to $N(N-1)>2$, a  condition
holds for any $N\ge2$.

On the other hand, it is easy to see that $\|Dv\|_{L^2(Q)}$ is uniformly bounded by testing the system with $v$. Assuming that all weak solutions of \mref{fullsys0a} are uniformly bounded and satisfy \mref{contcond}, we can start our argument  with $\frac{N^2}{2(N+2)}<p\le 1$  to  obtain $2p\cg >N$ and $2p\le2$. Thus, \reflemm{unifp} applies here with $q_0>q=p\cg$. Such $p\in (\frac{N^2}{2(N+2)},1]$ exists  if and only if  $\frac{N^2}{N+2}<2$ which is equivalent to $(N-1)^2<5$ or $N\le 3$. The spectral gap condition holds for $q$ if $\nu_*>1-1/(p\cg)$. We can choose $p$ such that $p\cg\sim N/2$ if $\nu_*>1-2/N$.
\eproof

Proving that a bounded weak solution of a cross diffusion system is uniformly continuous is already a hard problem. By the definition of $X$, a fixed point of $\mL$ in $X$ is H\"older continuous but their continuity is not uniform in order \mref{contcond} is verified so that we can apply \reflemm{Dup}. For any $N$ we can choose $\cg,\ag$ such that $\cg\lesssim 1+2/N$ and $\ag\gtrsim 1/\cg$ so that if $v\in X$ then $Dv\in L^{2q_0}(\Og)$ with $2q_0=\frac{2}{\cg(1-\ag)}>N$. But $\|Dv\|_{L^{2q_0}(\Og)}$ is not uniformly bounded so that the condition \mref{contcond} is not verified uniformly.

However, when $N=2$ the assumptions of \refcoro{Nlessequal3} can be greatly relaxed. In fact, the condition $\nu_*>1-2/N$ is clear and we can drop the condition \mref{contcond} to have

\bcoro{Nlessequal2} If $N= 2$  and solutions of \mref{fullsys0a} in $X$ are uniformly bounded then  \mref{fullsys0a} has a unique strong solution on $\Og\times(0,T)$. \ecoro

\bproof 
Let  $W$ be a fixed point of $\mL$. Again, testing the system with $W$ we see that if $\|W\|_{L^\infty(Q)}\le M$ then
$$\itQ{Q}{|DW|^2}\le C(M).$$

By \reflemm{lsu} ($f=\Div(\hat{\mb}(W))+\mg(W)$, so that $\|f\|_{L^2(Q)}\le C(M)$), we get $$\itQ{Q}{|D^2W|^2}\le C(M).$$

From the Gagliardo-Nirenberg BMO inequality (via a fixed covering)  we have
$$\itQ{Q}{|DW|^4}\le C\sup\|W\|_{BMO(\Og)}^2\itQ{Q}{|D^2W|^2}+C\frac{\sup\|W\|^2_{BMO(\Og)}}{R_0^2}\itQ{Q}{|DW|^{2}}.$$
Since $\|W\|_{BMO(\Og)}\le C\|W\|_{L^\infty(\Og)}$, together with \mref{defX2} (for $q=1$) we have
$$\iidx{\Og}{|DW|^{2}}+\itQ{Q}{|D^2W|^2}\le C\itQ{Q}{(|DW|^4+|DW|^2)}\le C(M)$$
for some constant $C(M)$. Therefore,
$$\iidx{\Og}{|DW|^{2}}\le C(M).$$
Hence, if $\|W\|_{L^\infty(Q)}\le M$ uniformly then we also have the  above estimate is uniform. We will show below that this implies $\|W\|_{BMO(\Og_R)}$ can be very small for uniformly small $R$. Of course, as $DW\in L^4(Q)$ and $4>N=2$. \reflemm{unifp} applies to give the result.

We now show that the uniform bound of $\|DW\|_{L^2(\Og)}\le M$ gives $\|W\|_{BMO(\Og_R)}$ small for uniformly small $R$ among weak solutions of \mref{fullsys0a}.

By contradiction, there is a sequence of weak solutions  $\{W_n\}$ converges to $W$ in $L^2(\Og)$  with $\|DW\|_{L^2(\Og)}\le M$. Furthermore, there are $\eg_0>0$ and sequences of positive numbers $\{r_n\}$ and points $ \{x_n\}\subset \Og$ such that $r_n\to0$ and $\|W_n\|_{BMO(\Og\cap B_{r_n}(x_n))}>\eg_0$. As $\bar{\Og}$ is compact, we can assume that $x_n\to x_0$ for some $x_0\in\Og$ and the balls $B_{r_n}(x_n)$ are concentric.

For any $R>0$ as $W_n\to W$ in $L^2(\Og)$
$$\mitx{\Og_R}{|W_n-(W_n)_R|}\to \mitx{\Og_R}{|W-W_R|}.$$

Since $DW\in L^2(\Og)$, for any $\eg>0$ we have $\|DW\|_{L^2(\Og_R)}<\eg$ if $R$ small. By Poincar\'e's inequality, $N=2$, we have $\|W\|_{BMO(\Og_R)}$ small for uniformly small $R$. For $r_n<R$ we have $\|W_n\|_{BMO(\Og_{r_n})}\le \|W_n\|_{BMO(\Og_{R})}$ which can be very small by the above limit. We then obtain a contradiction. \eproof

In the sequel, we will consider a special class of cross diffusion systems with polynomial growth parameters. The boundedness assumption of weak solutions will be greatly relaxed.  Inspired by the (SKT) system, we suppose that there is $k>0$ such that  $\llg(W)\sim |W|^k+1$ and
\beqno{cond4a} |\mb(W)|\le C\llg(W),\; |\hat{\mb}(W)| \mbox{ and }\mg(W)|\le C\llg(W)|W|.\eeq

To begin, we have the following elementary lemma showing that the energy is uniformly bounded under a very weak condition.

\blemm{Wboundlem} Assume the condition \mref{cond4a} and $N$ is any number.  If there is $M$ such that for any weak solution of \mref{fullsys0a} $$\itQ{Q}{|W|}\le M.$$ Then there is a constant $C(M)$ such that
$$\itQ{Q}{|W|^{k+2}}\le C(M),$$
and $$\itQ{Q}{|DW|^{2}}\le C(M).$$
\elemm

\bproof Test the system \mref{fullsys0a} with $W$ and use the growth condition \mref{cond4a} to easily get
$$\itQ{Q}{(1+|W|^k)|DW|^2}
\le C\itQ{Q}{|W|^{k+2}}+C.$$ This is to say for $V=|W|^{k/2+1}$
$$\itQ{Q}{|DV|^2}
\le C\itQ{Q}{|V|^{2}}+C.$$

By the interpolation inequality
$$\iidx{\Og}{|V|^2}\le \eg\iidx{\Og}{|DV|^2}+C(\eg,\bg)\left(\iidx{\Og}{|V|^\bg}\right)^\frac{2}{\bg}$$ we can choose $\eg,\bg$ sufficiently small to  have that
$$\itQ{Q}{(1+|W|^k)|DW|^2}
\le C\left(\itQ{Q}{|W|}\right)+C.$$
This proves the lemma. \eproof

For any $N$ we can choose $\cg,\ag$ such that $\cg\lesssim 1+2/N$ and $\ag\gtrsim 1/\cg$ so that if $v\in X$ then $Dv\in L^{2q_0\cg}(Q)$ with $2q_0\cg=\frac{2}{1-\ag}>N$. But the continuity condition \mref{contcond} is not verified easily even when $N=2$ and we have to assume it.

However, the essential role of \mref{contcond} is to absorb the integral of $|Dv|^{2p+2}$ on the right of \mref{keyassDv1} to the left in the proof of \reflemm{Dup}. Without \mref{contcond}, we have the following special cases in which the same purpose is served to obtain the same result where the boundedness assumption of solutions is replaced by that of their $BMO$ norms.

\bcoro{avsmall} Suppose $N\le3$ and $\nu_*>1/6$ if $N=3$. Suppose that the solutions of \mref{fullsys0a} in $X$ have uniform bounded BMO norms. If $\sup_{v\in\RR^m}\frac{|\ma_v|^2}{\llg(v)}$ is sufficiently small ($k\le 2$) then there exits strong  solutions in $\Og\times(0,T)$. 
\ecoro

\bproof We can choose $\cg$ near $1+2/N$ and $\ag$ near 1 in the definition of $X$. We then have $2q_0\cg=\frac{2}{1-\ag}\ge 4$. Thus $q_0\ge 2/\cg\sim 2N/(N+2)$. The spectral gap condition $\nu_*>1-1/q_0$ becomes $\nu_*>(N-2)/(2N)$ which is void if $N=2$ and $\nu_*>1/6$ if $N=3$.

Because $\sup_{(0,T)}\|v(\cdot,t)\|_{BMO(\Og)}<\infty$ we have $v\in L^r(Q)$ for all $r>1$. Therefore, the matrix parameters $\ma,\hat{\mb}, \mb$ and $\mg$ (and their derivatives in $v$) of the system are uniformly bounded in $L^q(Q)$ for any given $q>1$.

As in \reflemm{Dup}, $v$ weakly solves
$$(Dv)_t-\Div(\ma D(D v)+(\ma_vDv)Dv)=\Div(\hat{\mb}_vDv)+(\mb_vDv)Dv+\mb D(Dv)+\mg_v  Dv.$$
We formally test this system with $Dv\fg^2$ (see \reflemm{Dup}),  $\fg$ being a positive $C^1$ function. Thus, if $v\in X$ then $Dv\in L^4(Q)$ ($2q_0\cg\ge4$ in the definition of $X$) so that all the terms in \mref{keyassDv1} are finite (but they are {\em not uniformly} bounded) to get for any fixed point $v\in X$
$$\itQ{Q}{\llg(v)|D^2v|^2\fg^2}\le \itQ{Q}{|\ma_v(v)||Dv|^{2}|D^2v|\fg^2}+ C\itQ{Q}{|Dv|^{2}(\fg^2+|D\fg|^2)}.$$

By Young's inequality for some appropriate $q>1$
$$\itQ{Q}{\llg(v)|D^2v|^2\fg^2}\le \itQ{Q}{\frac{|\ma_v(v)|^2}{\llg(v)}|Dv|^{4}\fg^2}+\eg\itQ{Q}{|Dv|^4\fg^2}+C(\eg)\itQ{Q}{(\fg^q+|D\fg|^q)}.$$

For $C_*=\sup_{v\in\RR^m}\frac{|\ma_v|^2}{\llg(v)}+\eg$ and $\fg$ is a cutoff function for any ball $B_{R_0},B_{2R_0}$
$$\itQ{Q}{\llg(v)|D^2v|^2\fg^2}\le C_*\itQ{Q}{|Dv|^{4}\fg^2}+
C(\eg,R_0^{-1}).$$

Applying the Gagliardo-Nirenberg BMO inequality as in \reflemm{Dup} to the first term on the right,  we obtain for $C_{**}=C_*\|v\|_{BMO(\Og_{R_0})}^2$ and any  $R_0>0$
\beqno{keyassDv1b}\itQ{Q}{\llg(v)|D^2v|^2\fg^2}\le C_{**}\itQ{Q}{|D^2v|^{2}\fg^2}+C(C_{**},R_0^{-1})\itQ{Q}{|Dv|^{2}\fg^2}.\eeq

Thus, as  $\|v\|_{BMO(\Og)}$ is uniformly bounded, if  $\sup_{v\in\RR^m}\frac{|\ma_v|^2}{\llg(v)}$  is sufficiently small then we can absorb the first term on the right into that on the left to see that $\|D^2v\|_{L^2(\Og)}$ is bounded (a fact we cannot obtain by \reflemm{lsu} because $v$ is not assumed to be bounded) via a fixed covering (we fix  $R_0>0$). Arguing and iterating as in \reflemm{Dup} and \reflemm{unifp} we obtain a uniform bound for fixed points of $\mL$. In fact, we see now that $\|Dv\|_{L^4(\Og)}$ is uniformly bounded (by the Galiardo-Nirenberg inequality) and because $4>N$ if $N\le3$, the proof is complete. \eproof

As long as $C_{**}=\|v\|_{BMO(\Og_R)}\sup_{v\in\RR^m}\frac{|\ma_v|^2}{\llg(v)}$ in the proof of \refcoro{avsmall} is small for a fixed small $R>0$ then the proof can go on and we can conclude that

\bcoro{avsmallN2} Suppose $N=2$. Let $v$ be a strong solution in $\Og\times(0,T_0)$. Suppose that $\sup_{(0,T_0)}\|Dv\|_{L^2(\Og)}$ is bounded. If $\sup_{v\in\RR^m}\frac{|\ma_v|^2}{\llg(v)}$ is bounded then $v$ exists on $\Og\times(0,T_0)$. 
\ecoro

We then have the following result on the SKT system \mref{skt0a}.
\beqno{skt} \left\{\barr{ll}u_t= \Delta (P(u)) +\mg(u)&\mbox{in $\Og\times(0,T_0)$,}\\\mbox{Dirichlet or Neumann boundary conditions}&\mbox{on $\partial\Og\times(0,T_0)$,}\\u=u_0&\mbox{on $\Og$}.\earr\right.\eeq

\btheo{sktN2} Let $N=2$ and $u_0\in W^{1,2p}$ for some $p>1$. 
The problem \mref{skt} has a unique strong solution in $\Og\times (0,T_0)$ for any $T_0>0$.
\etheo

\bproof  Testing the system \mref{skt} with $u$, we easily obtain
\beqno{skt0}\|u\|_{L^2(\Og)},\; \|Du\|_{L^2(\Og\times(0,T_0))}\le C(\|u_0\|_{L^2(\Og)}).\eeq

We multiply the $i^{th}$ equation of \mref{skt} with $(P_i(u))_t$ and add the results to get
$$\iidx{\Og}{\myprod{u_t,P_u u_t}}+\frac12\frac{d}{dt} \iidx{\Og}{|DP|^2}=\iidx{\Og}{\myprod{\mg,P_u u_t}}.$$

Since $P_u$ is positive definite  (thanks to the ellipticity condition of $\ma$) and the components of $P,\mg$ are quadratic in $u$ we have
\beqno{skt1a}\frac{d}{dt} \iidx{\Og}{|DP|^2}+\iidx{\Og}{\llg(u)|u_t|^2}\le C\iidx{\Og}{(1+|u|^3)|u_t|}.\eeq
By Young's inequality, because $\llg(u)\ge \llg_0$ we derive 
\beqno{skt1}\frac{d}{dt} \iidx{\Og}{|DP|^2}\le C\iidx{\Og}{(1+|u|^6)}.\eeq

By the interpolation Gagliardo-Nirenberg inequality  $\|u\|_{L^6(\Og)}^6\le C\|u\|_{H^1(\Og)}^4\|u\|_{L^2(\Og)}^2$ when $N=2$.  From the ellipticity condition $\myprod{P_u x,x}\ge \llg_0\|x\|^2$ it easy to see that the matrix norm $\|P_u^{-1}\|\le 1/\llg_0$ so that, as $Du=P_u^{-1}DP(u)$, we have $|Du|\le \llg_0^{-1}|DP|$. Therefore, thanks to the estimate \mref{skt0} for $\|u\|_{L^2(\Og)}$
$$\|u\|_{L^6(\Og)}^6\le C(1+\|Du\|_{L^2(\Og)}^4)\le C\|Du\|_{L^2(\Og)}^2\iidx{\Og}{|DP|^2}+C,$$
where $C$ is a constant depends on $\|u_0\|_{L^2(\Og)}$.

Thus, for $y(t)=\iidx{\Og\times\{t\}}{|DP|^2}$ and $q(t)=C\|Du\|_{L^2(\Og)}^2$ we get from \mref{skt1}
$$y'(t)\le q(t)y(t)+c_0.$$

Here, $c_0$ is a constant depends on $\|u_0\|_{L^2(\Og)}$. This implies $$y(t)\le e^{\int_0^t q(s)ds}\left(y(0)+c_0\int_0^te^{-\int_0^s q(\tau)d\tau}ds\right).$$

By \mref{skt0} $\int_0^tq(s)ds \le C\itQ{\Og\times(0,T_0)}{|Du|^2}\le C(\|u_0\|_{L^2(\Og)})$, we see that for all $t\in(0,T_0)$
$$\iidx{\Og\times\{t\}}{|Du|^2} \le C\iidx{\Og\times\{t\}}{|DP|^2}\le C(T,\|u_0\|_{W^{1,2}(\Og)}).$$

Thus, $\|Du\|_{L^2(\Og)}$ is bounded on $(0,T_0)$ so that $\|u\|_{BMO(\Og_R)}$ is small if $R$ is uniformly small. By \reftheo{avsmallN2} the strong solution $u$ exists on $(0,T_0)$. \eproof

\bibliographystyle{plain}

\begin{thebibliography}{10}



\bibitem{Am2} H. Amann. \newblock{Dynamic theory of quasilinear parabolic systems III. Global existence,} {\em  Math Z.} 202 (1989), pp. 219–-250.

\bibitem{brezis} H. Brezis, \newblock{\em Functional Analysis, Sobolev Spaces and Partial Differential Equatiosn,} Universitext, Springer, New York, 2010. 

\bibitem{BC} O H. Brezis and M. Crandall, \newblock Uniqueness of solutions of the initial value problem for $u_t - \Delta \fg (u) = 0$, {\em J. Math. Pures Appl.} 58 (1979), p. 153--163.



\bibitem{AF} A. Friedman, \newblock{\em Partial Differential Equations,} New York, 1969.


\bibitem{JS} O. John and J. Stara. \newblock On the regularity of weak solutions to parabolic systems in two spatial dimensions. \newblock{\em Comm. P.D.E.}, 27(1998), pp. 1159–1170.

\bibitem{kuf1} K. H. W. K\"ufner,
\newblock{ Invariant regions for quasilinear reaction-diffusion systems and applications to a two population model}, {\em NoDEA}, 3(1996), 421--444.


\bibitem{GiaS} M. Giaquinta and M. Struwe. \newblock{ On the partial regularity of weak solutions of nonlinear parabolic
systems}. {\em  Math. Z.}, 179(1982),  437--451.





\bibitem{LSU} O. A Ladyzhenskaya, V. A. Solonnikov and N. N. Uraltseva], \newblock{\em Linear and Quasi-linear Equations of Parabolic Type,} Translations of Mathematical Monographs, AMS, 1968.




\bibitem{letrans} D. Le. \newblock Regularity of BMO weak solutions to nonlinear parabolic systems via homotopy. \newblock{Trans. Amer. Math. Soc}. 365 (2013), no. 5, 2723--2753.



\bibitem{LDg} D. Le. \newblock {Global Existence for Large Cross Diffusion Systems on Planar Domains}. {\em submitted}.


 

\bibitem{dleANS} D. Le. \newblock {Weighted Gagliardo-Nirenberg Inequalities Involving BMO Norms and Solvability of Strongly Coupled Parabolic Systems}. {\em Adv. Nonlinear Stud.} Vol. 16, No. 1(2016), 125--146.

\bibitem{dlebook} D. Le, \newblock{\em Strongly Coupled Parabolic and Elliptic Systems: Existence and Regularity of Strong/Weak Solutions.} De Gruyter, 2018.


\bibitem{LN} D. Le and V. Nguyen. \newblock{Global and blow up solutions to cross diffusion systems on 3D domains}, {\em Proc. AMS.} Vol. 144, No.11 (2016), 4845--4859.

\bibitem{dleJMAA} D. Le, \newblock{On the global existence of a generalized Shigesada-Kawasaki-Teramoto system,} {\em J. Math. Anal. App.} to appear.

\bibitem{lduni} D. Le \newblock{Uniqueness and Regularity of Unbounded Weak Solutions to a Class of Cross Diffusion Systems.} \newblock { 	arXiv:1906.03456}

\bibitem{LM} T. Lepoutre and A. Moussa. \newblock{Entropic structure and duality for multiple species cross-diffusion systems}, {\em Nonlinear Analysis} Vol. 159, (2017), 298--315.



\bibitem{NSver} J. Necas and V. Sverak. \newblock {On regularity of solutions of nonlinear parabolic systems.} {\em Ann. Scuola Norm. Sup. Pisa Cl. Sci.} (4), 18 (1), 1-11 (1991).



\bibitem{Red} R. Redlinger. \newblock{Existence of the global attractor for a strongly coupled parabolic system arising
in population dynamics}. {\em J. Diff. Eqns.}, 118(1995), 219--252.


\bibitem{SKT} N. Shigesada, K. Kawasaki  and E.~Teramoto. \emph{ 
Spatial segregation of interacting species}.  J. Theor. Biol., 
79(1979), 83-- 99.



\bibitem{yag} A. Yagi. \newblock{Global solution to some quasilinear parabolic systems in population dynamics}. {\em Nonlin. Anal}. 21
(1993), 603-630.
\end{thebibliography}

\end{document}